\documentclass[10pt]{article}
\usepackage{amsmath,amssymb}
\usepackage[margin=1in]{geometry}
\usepackage{booktabs}
\usepackage[colorlinks=true,linkcolor=blue,citecolor=blue,urlcolor=blue]{hyperref}

\title{Odd square roots and the sum of element orders on $S_n$ and $A_n$}
\author{Manoj Kumar Singh\\
\small Department of Mathematics and Statistics, St.~Xavier's College, Ranchi,\\
\small Jharkhand 834001, India\\
\small \texttt{manojkumarsingh@sxcran.org}}
\date{}

\begin{document}

\maketitle

\begin{abstract}
For a finite group $G$ let $\psi(G)$ be the sum of the orders of its elements. The quotient
$S_n/A_n$ has order two. Hence $\psi(S_n/A_n)=3$. The quantity to be compared with
$\psi(S_n)$ is therefore $3\psi(A_n)$. Computation indicates that $3\psi(A_n)>\psi(S_n)$
for every $n\ge3$. We reduce this inequality to a question on the parity of square roots.
Let $r(\beta)$ be the number of odd permutations $\sigma$ with $\sigma^{2}=\beta$. We show
that the sum of the orders of the odd permutations equals
$2\sum_{\beta\in A_n}r(\beta)o(\beta)$. We also show that $\sum_{\beta\in A_n}r(\beta)$
equals $|A_n|$. The inequality then says that $r$ and the order function are negatively
correlated on $A_n$. It says equally that the average order of an odd permutation is less
than twice the average order of an even one. We then rule out two natural approaches. No
injection from the odd permutations into $A_n$ can halve the order at every point. The
smallest failure occurs at $n=12$. A threshold argument through Landau's function fails as
well. Computations up to $n=60$ are reported.
\end{abstract}



\section{Introduction}\label{sec:intro}

Let $G$ be a finite group. For $g\in G$ let $o(g)$ be the order of $g$. Amiri, Jafarian
Amiri and Isaacs \cite{AAI2009} introduced the function
\begin{equation}\label{eq:psidef}
  \psi(G)=\sum_{g\in G}o(g).
\end{equation}
They proved that the cyclic group $\mathbb{Z}_n$ carries the largest value of $\psi$ among
all groups of order $n$. Later work used $\psi$ to detect nilpotency, solvability and
related structural properties \cite{AA2011,BaniasadKhosravi2018,Herzog2018,Tarnauceanu2020}.
Explicit values are known for several families. T\u{a}rn\u{a}uceanu and Fodor \cite{TF2014}
settled the abelian case. Chew, Chin and Lim \cite{CCL2017} gave a recursive formula for the
same class. Lazorec \cite{Lazorec2024} treated the extraspecial and the generalised
extraspecial $p$-groups. Banisadr, Refaghat, Marefat and Vakili \cite{BRMV2025} treated the
dihedral, the dicyclic, the Frobenius and the two dimensional linear groups. A survey is
available in \cite{Herzog2023}.

Two recent papers stand close to the present setting. M.~Amiri \cite{Amiri2025} proved that
two finite abelian groups of the same order are isomorphic exactly when their values of
$\psi$ agree. T\u{a}rn\u{a}uceanu \cite{Tarnauceanu2026} proved the two sided estimate
\begin{equation}\label{eq:tarn}
  \psi(H)+|H|\bigl(\psi(G/H)-1\bigr)\ \le\ \psi(G)\ \le\ \psi(H)+|H|^{2}\bigl(\psi(G/H)-1\bigr)
\end{equation}
for a proper normal subgroup $H$ of $G$. He also described the cases of equality.

The pair $(S_n,A_n)$ is a natural test case for a comparison of this kind. The quotient
$S_n/A_n$ is cyclic of order two. Hence $\psi(S_n/A_n)=3$. The quantity to be set against
$\psi(S_n)$ is $3\psi(A_n)$. Direct computation gives $\psi(S_3)/\psi(A_3)=13/7$ and
$\psi(S_4)/\psi(A_4)=67/31$. At $n=13$ the ratio equals
$118551513523/46287964867=2.56117\ldots$. In the whole computable range the ratio stays well
below $3$. This leads to the following statement.

\medskip\noindent
\textbf{Conjecture 1.1.} \emph{For every integer $n\ge3$,}
\begin{equation}\label{eq:main}
  3\,\psi(A_n)>\psi(S_n).
\end{equation}

\medskip
Our aim here is of two kinds. First we settle what \eqref{eq:main} really asserts. Then we
remove two approaches to it from consideration. Section~\ref{sec:reduction} reduces the
inequality to a statement on the parity of square roots. It then turns that statement into a
comparison of the order distributions of the two cosets of $A_n$ in $S_n$.
Section~\ref{sec:obstructions} carries the two negative results. An injection from the odd
permutations into $A_n$ that at least halves the order does not exist in general. The
smallest failure occurs at $n=12$. A threshold argument built on Landau's function has
nothing to work with. Section~\ref{sec:data} reports the computations. It also proposes a
sharper form of the conjecture.

Throughout, $n\ge2$ is an integer. Permutations act on $\{1,2,\dots,n\}$. Cycle types count
fixed points as cycles of length $1$. We write $\psi_{\mathrm{odd}}$ for the sum
$\sum_{\sigma\in S_n\setminus A_n}o(\sigma)$. Then
$\psi(S_n)=\psi(A_n)+\psi_{\mathrm{odd}}$. A cycle type of $n$ is a sequence
$(a_1,a_2,\dots)$ of non-negative integers with $\sum_L L\,a_L=n$, where $a_L$ counts the
cycles of length $L$. The conjugacy class of $S_n$ with this type has
\begin{equation}\label{eq:classsize}
  \frac{n!}{z},\qquad z=\prod_{L\ge1}L^{a_L}a_L!,
\end{equation}
elements. Every element of the class has order $\operatorname{lcm}\{L:a_L>0\}$. The elements
are even exactly when $n-\sum_L a_L$ is even.

\section{The reduction}\label{sec:reduction}

Squaring acts on cycles in two ways. A cycle of odd length $L$ squares to a cycle of the
same length $L$ on the same support. In the reverse direction an $L$-cycle with $L$ odd has
exactly one square root that is again an $L$-cycle on that support. It is the
$\tfrac{L+1}{2}$-th power of the cycle. A cycle of even length $2L$ squares to a pair of
disjoint $L$-cycles. In the reverse direction a pair of disjoint $L$-cycles arises in exactly
$L$ ways as the square of a $2L$-cycle on their joint support.

These two facts give the enumeration below. It is the cycle type form of the generating
function of Glebsky, Lic\'on and Rivera \cite{GLR2023}.

\medskip\noindent
\textbf{Proposition 2.1.} \emph{Let $\beta\in S_n$ have $a_L$ cycles of length $L$ for each
$L\ge1$. Then $\beta$ has a square root in $S_n$ if and only if $a_L$ is even for every even
$L$. In that case the square roots of $\beta$ correspond to the families
$(m_L)_{L\ge1}$ of non-negative integers with}
\begin{equation}\label{eq:mrange}
  m_L=\tfrac{1}{2}a_L\ \text{ for $L$ even},\qquad 0\le m_L\le\bigl\lfloor\tfrac{1}{2}a_L\bigr\rfloor
  \ \text{ for $L$ odd}.
\end{equation}
\emph{The number of square roots attached to a given family is}
\begin{equation}\label{eq:rootcount}
  \prod_{L\ge1}\frac{a_L!}{m_L!\,2^{m_L}\,(a_L-2m_L)!}\;L^{m_L}.
\end{equation}
\emph{Such a square root has $\sum_L(a_L-m_L)$ cycles. It is an odd permutation exactly when
$n-\sum_L(a_L-m_L)$ is odd.}

\medskip\noindent
\textit{Proof.} Let $\sigma^{2}=\beta$. Each cycle of $\sigma$ of even length $2L$ gives two
$L$-cycles of $\beta$. Each cycle of $\sigma$ of odd length $L$ gives one $L$-cycle of
$\beta$. Distinct cycles of $\sigma$ give cycles with disjoint supports. Every cycle of
$\beta$ arises in this way.

Read the last paragraph backwards. Then $\sigma$ is fixed by the following data. For each
$L$ we split the $a_L$ cycles of length $L$ into $m_L$ unordered pairs to be merged into
cycles of length $2L$ together with $a_L-2m_L$ cycles left alone. For each pair we choose a
merging. A cycle of $\beta$ left alone must have odd length. The reason is that a cycle of
even length is not the square of a cycle of the same length. So every cycle of even length
is paired. This forces $a_L$ to be even for $L$ even. It gives \eqref{eq:mrange}.

Now we count. The number of ways of choosing $m_L$ disjoint unordered pairs from $a_L$
objects is $a_L!/(m_L!\,2^{m_L}(a_L-2m_L)!)$. Each pair can be merged in exactly $L$ ways.
Each unpaired cycle of odd length has exactly one admissible root. Multiplying the three
counts gives \eqref{eq:rootcount}. Finally $\sigma$ has $m_L$ cycles of length $2L$ together
with $a_L-2m_L$ cycles of length $L$ for each $L$. So it has
$\sum_L\bigl(m_L+a_L-2m_L\bigr)=\sum_L(a_L-m_L)$ cycles. Its sign is
$(-1)^{n-\sum_L(a_L-m_L)}$. $\square$

\medskip
For $\beta\in A_n$ let $r(\beta)$ be the number of odd $\sigma\in S_n$ with
$\sigma^{2}=\beta$. Proposition 2.1 shows that $r(\beta)$ depends only on the cycle type of
$\beta$. To obtain it we sum \eqref{eq:rootcount} over the families $(m_L)$ for which
$n-\sum_L(a_L-m_L)$ is odd.

\medskip\noindent
\textbf{Theorem 2.2.} \emph{Let $n\ge2$. Then}
\begin{equation}\label{eq:identity}
  \psi_{\mathrm{odd}}=2\sum_{\beta\in A_n}r(\beta)\,o(\beta)
  \qquad\text{and}\qquad
  \sum_{\beta\in A_n}r(\beta)=|A_n|.
\end{equation}
\emph{Hence \eqref{eq:main} holds for a given $n$ if and only if}
\begin{equation}\label{eq:corr}
  \sum_{\beta\in A_n}\bigl(r(\beta)-1\bigr)o(\beta)<0.
\end{equation}

\medskip\noindent
\textit{Proof.} Let $\sigma$ be an odd permutation. Then
$1=\operatorname{sgn}(\sigma^{o(\sigma)})=\operatorname{sgn}(\sigma)^{o(\sigma)}
=(-1)^{o(\sigma)}$. So $o(\sigma)$ is even. Therefore $\sigma^{2}$ lies in $A_n$ and
$o(\sigma^{2})=o(\sigma)/2$. Group the odd permutations by the value of $\sigma^{2}$. This
gives
\[
  \psi_{\mathrm{odd}}=\sum_{\sigma\ \text{odd}}o(\sigma)
  =\sum_{\sigma\ \text{odd}}2\,o(\sigma^{2})
  =2\sum_{\beta\in A_n}r(\beta)\,o(\beta),
\]
which is the first identity. The fibres of the map $\sigma\mapsto\sigma^{2}$ on the odd
permutations partition a set of size $|S_n\setminus A_n|=|A_n|$. This gives the second
identity.

We have $\psi(S_n)=\psi(A_n)+\psi_{\mathrm{odd}}$. So the inequality $3\psi(A_n)>\psi(S_n)$
is the same as $\psi_{\mathrm{odd}}<2\psi(A_n)$. By the first identity this reads
$\sum_{\beta}r(\beta)o(\beta)<\sum_{\beta}o(\beta)$. That is \eqref{eq:corr}. $\square$

\medskip
The second identity in \eqref{eq:identity} says that $r$ has average value $1$ on $A_n$.
Condition \eqref{eq:corr} then says that $r$ and $o$ are negatively correlated. In words,
the elements of large order must carry fewer odd square roots than the average.

The constant $3$ in \eqref{eq:main} is not a generous one. Suppose that $r(\beta)=1$ for
every $\beta\in A_n$. Then \eqref{eq:identity} gives $\psi_{\mathrm{odd}}=2\psi(A_n)$. So
$\psi(S_n)=3\psi(A_n)$ exactly. The average value of $r$ on $A_n$ equals $1$ in every case.
So the constant $3$ is the value produced by a complete lack of correlation between $r$ and
$o$. The remark is not idle. At $n=2$ the group $A_2$ is trivial. There $r$ takes the value
$1$ at the identity. Accordingly $\psi(S_2)=3=3\psi(A_2)$. This equality is the reason for
the range $n\ge3$ in Conjecture 1.1. The conjecture therefore carries no slack in its
formulation. The whole of the observed gap comes from a genuine negative correlation.

A cheap bound is not available here. From \eqref{eq:identity} one gets
$\psi_{\mathrm{odd}}\le2\bigl(\max_{\beta}r(\beta)\bigr)\psi(A_n)$. This settles
\eqref{eq:main} only when $\max_\beta r(\beta)\le1$. The maximum is in fact large. On the
class of cycle type $(8,8)$ in $A_{16}$ the value of $r$ equals $8$. In general $r$ grows
with the repeated cycle length.

The correlation can be restated through the level sets of the order function. For $t\ge1$
put
\[
  N_{\mathrm{odd}}(t)=\#\{\sigma\in S_n\setminus A_n:\ o(\sigma)\ge t\},\qquad
  N_{\mathrm{even}}(t)=\#\{\beta\in A_n:\ o(\beta)\ge t\}.
\]
Let $R(s)$ be the sum of $r(\beta)$ over the $\beta\in A_n$ with $o(\beta)\ge s$.

\medskip\noindent
\textbf{Proposition 2.3.} \emph{For every $t\ge1$ we have
$N_{\mathrm{odd}}(t)=R\bigl(\lceil t/2\rceil\bigr)$. Hence Conjecture 1.1 holds for a given
$n\ge3$ if and only if}
\begin{equation}\label{eq:levelsets}
  \sum_{s\ge1}\Bigl(R(s)-N_{\mathrm{even}}(s)\Bigr)<0 .
\end{equation}

\medskip\noindent
\textit{Proof.} An odd $\sigma$ satisfies $o(\sigma)=2o(\sigma^{2})$ with $\sigma^{2}$ in
$A_n$. This was shown in the proof of Theorem 2.2. So $o(\sigma)\ge t$ holds if and only if
$o(\sigma^{2})\ge t/2$ holds. Orders are integers. So the last condition is
$o(\sigma^{2})\ge\lceil t/2\rceil$. Grouping by the value of $\sigma^{2}$ gives the first
assertion.

For a non-negative integer valued function $f$ on a finite set we have
$\sum_{x}f(x)=\sum_{t\ge1}\#\{x:f(x)\ge t\}$. Applied to $o$ this gives
\[
  \psi_{\mathrm{odd}}=\sum_{t\ge1}N_{\mathrm{odd}}(t)=\sum_{t\ge1}R\bigl(\lceil t/2\rceil\bigr)
  =2\sum_{s\ge1}R(s).
\]
The last step holds because each $s$ occurs for the two values $t=2s-1$ and $t=2s$. In the
same way $\psi(A_n)=\sum_{s\ge1}N_{\mathrm{even}}(s)$. So the inequality
$\psi_{\mathrm{odd}}<2\psi(A_n)$ is exactly \eqref{eq:levelsets}. $\square$

\medskip
Divide \eqref{eq:corr} by $|A_n|=|S_n\setminus A_n|=n!/2$. The constant $3$ then disappears
from the statement.

\medskip\noindent
\textbf{Corollary 2.4.} \emph{Conjecture 1.1 is equivalent to the following assertion. For
every $n\ge3$ the average order of an odd permutation of $\{1,\dots,n\}$ is smaller than
twice the average order of an even permutation.}

\medskip
The termwise form of \eqref{eq:levelsets} is the condition $R(s)\le N_{\mathrm{even}}(s)$
for all $s$. It is strictly stronger. It says that $r$ has average at most $1$ on every
upper level set of $o$. Condition \eqref{eq:corr} asks for this only in a weighted sense.
The stronger form is false. We show this next.

\section{Two obstructions}\label{sec:obstructions}

We first record the classes on which $r$ vanishes. These classes carry the whole surplus in
\eqref{eq:corr}.

\medskip\noindent
\textbf{Lemma 3.1.} \emph{Let $\beta\in A_n$ have pairwise distinct cycle lengths. Fixed
points are counted here as cycles of length $1$. Then $r(\beta)=0$.}

\medskip\noindent
\textit{Proof.} Suppose $\sigma$ is odd with $\sigma^{2}=\beta$. A permutation with $a_j$
cycles of length $j$ has sign $(-1)^{\sum_j a_j(j-1)}$. So an odd permutation has an odd
number of cycles of even length. In particular $\sigma$ has a cycle of even length $2L$ for
some $L\ge1$. Squaring carries that cycle to two disjoint cycles of length $L$. Every other
cycle of $\sigma$ gives cycles supported inside itself. So $\beta$ has at least two cycles
of length $L$. This contradicts the hypothesis. $\square$

\medskip
The same conclusion follows from Proposition 2.1. If all $a_L\le1$ then the only admissible
family is $m_L=0$ for every $L$. The root attached to it has $\sum_L a_L$ cycles. So it has
the same sign as $\beta$.

Lemma 3.1 needs the fixed points to be counted. Without them the conclusion fails. The
element $\beta=(1\,2\,3)$ of $A_5$ has a single non-trivial cycle. Still $(1\,3\,2)(4\,5)$
is an odd square root of it. The next lemma explains this behaviour.

\medskip\noindent
\textbf{Lemma 3.2.} \emph{Let $\beta\in A_n$ have at least two fixed points. Suppose $\beta$
has a square root in $S_n$. Then $\beta$ has square roots of both parities. In particular
$r(\beta)>0$.}

\medskip\noindent
\textit{Proof.} Let $\sigma\in S_n$ satisfy $\sigma^{2}=\beta$. Let $F$ be the set of fixed
points of $\beta$, so that $|F|\ge2$. Take $a$ with $\beta(a)=a$. Then
$\beta(\sigma(a))=\sigma(\beta(a))=\sigma(a)$. So $\sigma$ maps $F$ into itself. The square
$\sigma^{2}$ fixes $F$ pointwise. So the restriction $\sigma|_F$ is an involution.

Suppose first that $\sigma|_F$ is not the identity. Choose $x,y\in F$ with $\sigma(x)=y$
and $\sigma(y)=x$. Then $\{x,y\}$ is an orbit of $\sigma$. So $\sigma=(x\,y)\tau$ with
$\tau$ fixing $x$ and $y$ and disjoint from $(x\,y)$. The two factors commute. Hence
$\tau^{2}=\sigma^{2}=\beta$.

Suppose next that $\sigma|_F$ is the identity. Choose any two distinct $x,y\in F$. The
transposition $(x\,y)$ is disjoint from the support of $\sigma$. Hence
$\bigl((x\,y)\sigma\bigr)^{2}=(x\,y)^{2}\sigma^{2}=\beta$. In both cases we produced two
square roots of $\beta$ that differ by a transposition. They have opposite parity. $\square$

\medskip
The hypothesis on the existence of a square root cannot be dropped. Take an element of $A_8$
of cycle type $(4,2,1,1)$. It has two fixed points. It has no square root at all. A
$4$-cycle of $\beta$ can only come from an $8$-cycle of $\sigma$. Squaring an $8$-cycle
produces two $4$-cycles.

By Theorem 2.2 every $\beta$ as in Lemma 3.1 contributes exactly $-o(\beta)$ to the left
side of \eqref{eq:corr}. Write
\begin{equation}\label{eq:creditdebit}
  C_n=\sum_{\substack{\beta\in A_n\\ r(\beta)=0}}o(\beta),
  \qquad
  D_n=\sum_{\substack{\beta\in A_n\\ r(\beta)\ge1}}\bigl(r(\beta)-1\bigr)o(\beta).
\end{equation}
Then Conjecture 1.1 for a given $n$ is the inequality $D_n<C_n$. Computation shows that
$C_n$ is dominated by the classes of Lemma 3.1. It shows that $D_n$ is dominated by the
classes with exactly one repeated cycle length. Take $n=13$ as an example. The largest
contributions to $C_{13}$ come from the cycle types $(13)$, $(7,5,1)$, $(10,2,1)$ together
with $(8,3,2)$. The largest contributions to $D_{13}$ come from $(6,6,1)$, $(5,5,3)$
together with $(5,4,4)$. Here $r$ takes the values $6$, $5$ and $4$ in that order.

\subsection*{Order halving injections}

We have $|S_n\setminus A_n|=|A_n|$. We also have $o(\sigma)=2\,o(\sigma^{2})$ for odd
$\sigma$. So one is tempted to prove Conjecture 1.1 by an injection. Let $\Psi$ map
$S_n\setminus A_n$ into $A_n$ with $o(\Psi(\sigma))\ge o(\sigma)/2$ for all $\sigma$, with
strict inequality for at least one $\sigma$. Summing over the odd permutations would give
$\psi_{\mathrm{odd}}<2\psi(A_n)$ at once. No such map exists in general.

\medskip\noindent
\textbf{Theorem 3.3.} \emph{Let $n\ge2$. An injection $\Psi:S_n\setminus A_n\to A_n$ with
$o(\Psi(\sigma))\ge o(\sigma)/2$ for every odd $\sigma$ exists if and only if}
\begin{equation}\label{eq:hall}
  N_{\mathrm{odd}}(t)\le N_{\mathrm{even}}\bigl(\lceil t/2\rceil\bigr)
  \qquad\text{for every }t\ge1 .
\end{equation}

\medskip\noindent
\textit{Proof.} Form the bipartite graph on $(S_n\setminus A_n)\cup A_n$. Join $\sigma$ to
$\beta$ when $o(\beta)\ge o(\sigma)/2$. The required map is a matching that saturates the
odd side. By Hall's theorem it exists if and only if $|N(S)|\ge|S|$ for every set $S$ of odd
permutations. Here $N(S)$ is the neighbourhood of $S$.

The neighbourhood of a single $\sigma$ is
$\{\beta\in A_n:o(\beta)\ge\lceil o(\sigma)/2\rceil\}$. It depends only on $o(\sigma)$.
These sets are totally ordered by inclusion. They decrease as $o(\sigma)$ grows. So for any
$S$ the set $N(S)$ agrees with the neighbourhood of an element of $S$ of maximal order. For
a fixed $N(S)$ the size $|S|$ is largest when $S$ consists of all odd permutations of order
at least $t=\max_{\sigma\in S}o(\sigma)$. So it is enough to test the sets
$S_t=\{\sigma\ \text{odd}:o(\sigma)\ge t\}$. For these the condition $|N(S_t)|\ge|S_t|$
reads exactly as \eqref{eq:hall}. $\square$

\medskip\noindent
\textbf{Theorem 3.4.} \emph{Condition \eqref{eq:hall} fails for $n=12$. Hence no injection
$\Psi:S_{12}\setminus A_{12}\to A_{12}$ satisfies $o(\Psi(\sigma))\ge o(\sigma)/2$ for all
$\sigma$.}

\medskip\noindent
\textit{Proof.} Take $t=24$. The odd classes of $S_{12}$ whose elements have order at least
$24$ are those of cycle types
\[
  (8,3,1),\quad(7,4,1),\quad(7,3,2),\quad(6,5,1),\quad(5,4,3),\quad(5,3,2,1,1).
\]
Their orders are $24$, $28$, $42$, $30$, $60$ and $30$. Formula \eqref{eq:classsize} gives
the sizes
\[
  19958400,\quad 17107200,\quad 11404800,\quad 15966720,\quad 7983360,\quad 7983360.
\]
These add up to $N_{\mathrm{odd}}(24)=80403840$. The classes of $A_{12}$ whose elements have
order at least $12$ are those of cycle types
\[
  (7,5),\ (5,4,2,1),\ (7,3,1,1),\ (6,4,1,1),\ (7,2,2,1),\ (5,3,3,1),\ (4,4,3,1),\
  (5,3,2,2),\ (4,3,3,2),\ (4,3,2,1,1,1),\ (5,3,1,1,1,1).
\]
Their sizes are
\[
  13685760,\ 11975040,\ 11404800,\ 9979200,\ 8553600,\ 5322240,\ 4989600,\ 3991680,\
  3326400,\ 3326400,\ 1330560.
\]
These add up to $N_{\mathrm{even}}(12)=77885280$. Now $80403840>77885280$. So
\eqref{eq:hall} fails at $t=24$. Theorem 3.3 gives the conclusion. $\square$

\medskip
The two totals can be checked against one another through Proposition 2.3. Summing
$r(\beta)$ over the eleven classes on the even side returns $80403840$. This agrees with
$N_{\mathrm{odd}}(24)=R(12)$.

Condition \eqref{eq:hall} fails at further places. For $n=12$ it fails also at $t=28$. For
$n=13$ it fails at $t=28$ and at $t=30$. For $n=18$ it fails at $t=180$. For $n=19$ it fails
at $t=420$. For $n=30$ it fails at one threshold. For every other $n$ with $2\le n\le40$ it
holds at all $t$. Conjecture 1.1 is still valid at each of these exceptional values. So
\eqref{eq:hall} is genuinely stronger than \eqref{eq:main}.

A relaxation of the injection to a transport plan does not help. A fractional matching that
saturates the odd side solves the linear programme attached to the same bipartite graph. The
bipartite matching polytope is integral. So a fractional matching exists exactly when an
integral one exists. Theorem 3.4 therefore rules out every argument that spreads the mass of
each odd permutation over even permutations of at least half its order without using the
numerical values of the orders any further.

\subsection*{Thresholds and Landau's function}

Let $g(n)$ be Landau's function, the largest order of an element of $S_n$. Let $g_2(n)$ be
the largest order of a permutation of $\{1,\dots,n\}$ that has at least two cycles of equal
length. Lemma 3.1 gives $r(\beta)=0$ whenever $o(\beta)>g_2(n)$. So one may try to split the
sum in \eqref{eq:levelsets} at $s=g_2(n)$. The range $s>g_2(n)$ has $R(s)=0$. It would then
serve as a reservoir of surplus. The reservoir turns out to be empty.

\medskip\noindent
\textbf{Proposition 3.5.} \emph{For every $n\ge3$ we have $g_2(n)\ge g(n-2)$. Moreover
$g_2(21)=g(21)=420$ and $g_2(22)=g(22)=420$.}

\medskip\noindent
\textit{Proof.} Let $\alpha$ be a permutation of $\{1,\dots,n-2\}$ of order $g(n-2)$. Regard
$\alpha$ as a permutation of $\{1,\dots,n\}$ that fixes the two remaining points. Its cycle
type now contains two cycles of length $1$. So it is counted by $g_2(n)$. Its order is still
$g(n-2)$. This gives the inequality. For the two numerical claims take the type
$(7,5,4,3,1,1)$ on $21$ points. It has order $420=g(21)$ together with a repeated cycle
length. Adjoining one more fixed point gives the claim for $22$. $\square$

\medskip
The values $g(n-2)$ and $g(n)$ agree for many $n$. They are always of the same order of
magnitude. So Proposition 3.5 leaves almost nothing above the threshold.
Table~\ref{tab:landau} records $g$ and $g_2$. It also records the share of $\psi(A_n)$ that
comes from the classes with $o(\beta)>g_2(n)$. That share is zero for most $n\ge14$ in the
computed range. From $n=15$ onwards it never exceeds $0.09$. A proof of Conjecture 1.1 must
therefore compare the two cosets in the range where $r$ is positive. It cannot isolate a
range where $r$ vanishes.

\begin{table}[ht]
\centering
\begin{tabular}{@{}rrrr@{}}
\toprule
$n$ & $g(n)$ & $g_2(n)$ & share of $\psi(A_n)$ above $g_2(n)$\\
\midrule
12 & 60 & 30 & 0.161\\
13 & 60 & 30 & 0.135\\
14 & 84 & 60 & 0.000\\
15 & 105 & 60 & 0.083\\
16 & 140 & 84 & 0.072\\
18 & 210 & 140 & 0.000\\
21 & 420 & 420 & 0.000\\
22 & 420 & 420 & 0.000\\
30 & 4620 & 2310 & 0.000\\
40 & 27720 & 16380 & 0.000\\
\bottomrule
\end{tabular}
\caption{Landau's function together with the largest order reached with a repeated cycle
length. The last column gives the share of $\psi(A_n)$ lying above that order.}
\label{tab:landau}
\end{table}

\section{Computations, a sharper conjecture and what remains}\label{sec:data}

The values below come from a sum over cycle types in exact integer arithmetic. The class
sizes are taken from \eqref{eq:classsize}. The counts $r$ are taken from Proposition 2.1.
The identity \eqref{eq:identity} was checked for $3\le n\le18$ by computing its two sides
from unrelated data. The left side came from the odd cycle types. The right side came from
the root counts. The two agree in every case.

\begin{table}[ht]
\centering
\begin{tabular}{@{}rrrr@{}}
\toprule
$n$ & $\psi(A_n)$ & $\psi(S_n)$ & $\psi(S_n)/\psi(A_n)$\\
\midrule
3 & 7 & 13 & 1.857\\
4 & 31 & 67 & 2.161\\
5 & 211 & 471 & 2.232\\
6 & 1411 & 3271 & 2.318\\
7 & 12601 & 31333 & 2.487\\
8 & 137047 & 299223 & 2.183\\
9 & 1516831 & 3291487 & 2.170\\
10 & 18111751 & 39020911 & 2.154\\
11 & 223179001 & 543960561 & 2.437\\
12 & 2973194071 & 7466726983 & 2.511\\
13 & 46287964867 & 118551513523 & 2.561\\
14 & 835826439631 & 1917378505407 & 2.294\\
15 & 15722804528341 & 32405299019941 & 2.061\\
16 & 292673102609791 & 608246253790591 & 2.078\\
\bottomrule
\end{tabular}
\caption{The ratio $\psi(S_n)/\psi(A_n)$ for small $n$.}
\label{tab:ratios}
\end{table}

Conjecture 1.1 was checked for all $n$ with $3\le n\le60$. On that range the ratio
$\psi(S_n)/\psi(A_n)$ never exceeds $2.56117\ldots$. That value is reached at $n=13$. The
ratio shows no tendency to grow. In the same range $\psi_{\mathrm{odd}}/\psi(A_n)$
oscillates within $[1.061,1.562]$ for $6\le n\le40$. So in the form of Corollary 2.4 there
is a factor of about $1.3$ to spare. The oscillation has an arithmetic source. It records
whether the cycle types of largest order at a given $n$ happen to be even or happen to be
odd. The quotient $D_n/C_n$ of \eqref{eq:creditdebit} stays within $[0.32,0.70]$ for
$3\le n\le28$. Its maximum again falls at $n=13$.

These computations suggest a sharper statement.

\medskip\noindent
\textbf{Conjecture 4.1.} \emph{For every $n\ge3$,}
\[
  \frac{\psi(S_n)}{\psi(A_n)}\ \le\ \frac{\psi(S_{13})}{\psi(A_{13})}
  =\frac{118551513523}{46287964867}=2.56117\ldots
\]
\emph{Equality holds only at $n=13$.}

\medskip
Conjecture 4.1 implies Conjecture 1.1. It also points to a change of the reference constant.
The natural comparison is not with the number $3$ coming from $\psi(S_n/A_n)$. It is with
the limiting behaviour of the ratio. The data place that limit near $2$. Such a value would
mean that the two cosets have asymptotically equal average order.

Theorem 2.2 with Corollary 2.4 reduces Conjecture 1.1 to a comparison of the order
distributions of the two cosets of $A_n$ in $S_n$. Theorem 3.4 shows that the comparison
cannot be made by a pointwise assignment that respects the factor $2$. Proposition 3.5 shows
that it cannot be made by isolating the elements of largest order. One route stays open. It
is analytic. One would prove that $\psi_{\mathrm{odd}}/\psi(A_n)$ is bounded away from $2$
using the distribution of $\log o(\sigma)$ on each coset separately. Results of the
Erd\H{o}s--Tur\'an type describe that distribution on $S_n$ as a whole. What is needed is a
form of them that is sensitive to the sign character. Such a form does not seem to be
available at present.

The computations raise two smaller questions. The first is to describe the integers $n$ for
which the strong condition \eqref{eq:hall} fails. Within $2\le n\le40$ these are $n=12$,
$13$, $18$, $19$ and $30$. No pattern is visible. The second is to decide Conjecture 4.1.

\subsection*{Acknowledgements}
The computations were carried out with the open source \textsf{Python} language. The author
thanks its developers.

\subsection*{Disclosure of interest}
The author reports there are no competing interests to declare.

\subsection*{Funding}
No funding was obtained for the reported work.

\subsection*{Data availability statement}
No data sets were generated or analysed during the current study. All computational results
are reproducible from the descriptions given in the paper.

\end{document}